\documentclass[11pt]{article}
\usepackage{amssymb,amsmath,latexsym}

\newtheorem{theo+}           {Theorem}      [section]
\newtheorem{prop+}  [theo+]  {Proposition}
\newtheorem{coro+}  [theo+]  {Corollary}
\newtheorem{lemm+}  [theo+]  {Lemma}
\newtheorem{exam+}  [theo+]  {Example}
\newtheorem{rema+}  [theo+]  {Remark}
\newtheorem{defi+}  [theo+]  {Definition}

\newtheorem{exam+s}  [theo+]  {Examples}
\newtheorem{rema+s}  [theo+]  {Remarks}
\newtheorem{hyp+}  [theo+]  {Hypotheses}
\newtheorem{cla+}  [theo+]  {Claim}

\newenvironment{proposition}{\begin{prop+}}{\end{prop+}}

\newenvironment{lemma}{\begin{lemm+}}{\end{lemm+}}

\def\bull{$\Box$}
\def\endmark{\hskip 2em\bull\par}
\newcommand{\qed}{\null\hspace*{\fill}\endmark}
\newcommand{\proof}{\noindent {\bf Proof} \hskip 0.4em}
\newcommand{\eproof}{\bigskip}

\newcommand{\ds}{\displaystyle}
\newcommand{\ts}{\textstyle}
\newcommand{\pa}{\partial}
\newcommand{\half}{{\textstyle\frac{1}{2}}}

\newcommand{\abs}[1]{\lvert#1\rvert} %modulus of a number
\newcommand{\inn}[1]{\langle#1\rangle} %inner product
\newcommand{\biginn}[1]{\big\langle#1\big\rangle} %\big version
\newcommand{\Biginn}[1]{\Big\langle#1\Big\rangle} %\Big version
\renewcommand{\vec}[1]{\boldsymbol{#1}} %for vectors (\vec usually
\newcommand{\RR}{{\mathbb R}}  %real numbers
\newcommand{\CC}{{\mathbb C}}  %complex numbers

\newcommand{\CP}{{\mathbb C}P} %complex projective space

\newcommand{\cc}{\mathbf c} %on derivs etc
\newcommand{\cC}{\mathbf C} %on inner products etc
\newcommand{\rR}{\mathbf R} %on diffls etc
\newcommand{\hh}{\text{Herm}}

\renewcommand{\phi}{\varphi}
\newcommand{\la}{\lambda}

\newcommand{\na}{\nabla}

\newcommand{\Ga}{\Gamma}

\newcommand{\om}{\omega}

\newcommand{\ov}{\overline} %often better than \bar
\DeclareMathOperator{\Tr}{Trace} %Trace

\newcommand{\dd}{{\mathrm d}}  %differential
\newcommand{\DD}{\nabla}  %covariant differential

\newcommand{\dt}{\ensuremath{\left.\frac{\dd}{\dd t}\right\vert_{0}}}
\newcommand{\pat}{\ensuremath{\left.\frac{\pa}{\pa t}\right\vert_{0}}}
\newcommand{\Dt}{\ensuremath{\left.\frac{\DD}{\pa t}\right\vert_{0}}}
\newcommand{\dcphi}{\ensuremath{\frac{\pa^{\cc}\phi}{\pa z}}}
\newcommand{\dphi}{\ensuremath{\frac{\pa\phi}{\pa z}}}
\newcommand{\dcbarphi}{\ensuremath{\frac{\pa^c\phi}{\pa\bar{z}}}}
\newcommand{\dbarphi}{\ensuremath{\frac{\pa\phi}{\pa\bar{z}}}}
\newcommand{\dbar}{\ensuremath{\frac{\pa}{\pa\bar{z}}}}
\newcommand{\Dz[1]}{\ensuremath{\frac{\DD^{#1}}{\pa z^{#1}}}}
\newcommand{\Dzbar[1]}{\ensuremath{\frac{\DD^{#1}}{\pa\bar{z}^{#1}}}}
\newcommand{\dcphit}{\ensuremath{\frac{\pa^{\cc}\phi_t}{\pa z}}}

\begin{document} 
\title{Jacobi fields along harmonic maps\thanks{Paper presented
at the Mathematical Society of Japan 9th International Research Institute: 
`Integrable Systems in Differential Geometry',
17-21 July 2000, University of Tokyo.}}

\author{John C.~Wood
 \\
{\small Department of Pure Mathematics, University of Leeds, 
Leeds LS2 9JT, G.B.} \\
{\small e-mail address: j.c.wood@leeds.ac.uk} }

\date{} 
\maketitle

\begin{abstract}
We show that Jacobi fields along harmonic maps between suitable spaces
preserve conformality, holomorphicity, real isotropy and complex isotropy
to first order;
this last being one of the key tools in the proof by Lemaire and the author
of integrability of Jacobi fields along harmonic maps from the $2$-sphere to
the complex projective plane.

\bigskip
\noindent\emph{Keywords:}  Harmonic map; minimal surface; Jacobi field

\noindent\emph{2000 MSC:} 58E20, 53C43, 53A10
\end{abstract}

\title{}

\section{Introduction}

A \emph{Jacobi field} (for the energy) along a harmonic map $\phi$ is
a vector field along
$\phi$ which is in the kernel of the second variation of the energy.
Equivalently, it
is tangent to a variation of $\phi$ for which the tension field remains
zero `to first order'.   We shall show that Jacobi fields preserve
several properties of a harmonic map to first order, namely
holomorphicity for maps between compact K\"ahler manifolds, conformality
of maps from the $2$-sphere, real isotropy of a harmonic map from the
$2$-sphere to a space form, and complex isotropy of a harmonic map from
the $2$-sphere to a complex space form.   The main idea of the proofs is
to define differentials depending on the map and the Jacobi field which
are shown to be holomorphic and so vanish; this generalizing the proof
of isotropy of a harmonic map from the $2$-sphere as
given, for example, in \cite{Wood-Horwood}.

Note that any harmonic map from the $2$-sphere is weakly conformal
and so is the same thing as
a \emph{minimal branched immersion} in the sense of \cite{G-O-R}.

Lastly, we mention that the Jacobi equation is preserved by
\emph{harmonic morphisms}, see \cite{MonWoo} for this result,
\cite{hamobib} for a bibliography and
\cite{Book} for an account of the theory of harmonic morphisms.

The author thanks the organizers of the MSJ 9th International Research
Institute, Tokyo 2000, for inviting him
to present this paper, and Luc Lemaire for some useful comments on it.

\section{Jacobi fields}

Throughout this paper $M = (M,g)$ and $N= (N,h)$ will denote smooth
(i.e. $C^{\infty}$) compact Riemannian manifolds without boundary and
$\phi:M \to N$ will denote a smooth map.
For any vector bundle $E \to M$, $\Ga(E)$ will denote the space of smooth
sections of $E$.   Particularly important is the \emph{pull-back bundle}
$\phi^{-1}TN \to M$; its smooth sections are called
\emph{vector fields along $\phi$}.

We recall some basic definitions, see \cite{Topics} for more details.
The \emph{energy} of $\phi$ is
defined by the integral
$$
E(\phi) = \half \int_M \abs{\dd\phi}^2 \,\om_g
$$
where $\abs{\dd\phi}$ denotes the Hilbert--Schmidt norm of the
differential of $\phi$
%given in local coordinates $(x^i)zz$ by
%$\abs{\dd\phi}^2 = g^{ij} h_{\al\be}
%\frac{\pa\phi^{\al}}{\pa x^i}\frac{\pa\phi^{\be}}{\pa x^j}$
and $\om_g$ denotes the volume form of the metric $g$.
By a \emph{(smooth) ($1$-parameter) variation} $\{\phi_t\}$ of $\phi$ we mean a smooth map
$M \times I \to N$, \ $(x,t) \mapsto \phi_t(x)$, where $I$ is an open
interval of the real line containing $0$,
such that $\phi_0 = \phi$.

Given a smooth variation $\{\phi_t\}$, we set
$$
v = \left.\frac{\pa\phi_t}{\pa t}\right\vert_{0}
$$
where $\pat$ denotes (partial) derivative with respect to $t$ at
$t=0$;
this defines a vector field along $\phi$
called the \emph{variation vector field of $\{\phi_t\}$}.
Then, for any smooth map $\phi:M \to N$ and any smooth variation
$\{\phi_t\}$ of it we have
\begin{equation} \label{1st-var1}
\dt E(\phi_t) = - \int_M
\inn{\tau(\phi), v} \,\om_g
\end{equation}
where $\tau(\phi)$ denotes the vector field along $\phi$ called the
\emph{tension field} of
$\phi$ given by $\tau(\phi) = \Tr\na\dd\phi$.

Here and in the sequel, $\langle \ , \ \rangle$ and $\nabla$
denote the inner product and connection induced on the relevant
bundle by the metrics and Levi-Civita connections on $M$ and $N$,
see, for example, \cite{Topics} for this formalism.

The formula \eqref{1st-var1} shows that the number
$\dt E(\phi_t)$ depends
only on the variation vector field $v$ of $\{\phi_t\}$.
Given any $v \in \Ga(\phi^{-1}TN)$, there are infinitely many smooth variations
$\{\phi_t\}$ \emph{tangent to} (i.e.\ \emph{with variation vector field}) $v$.
Writing the left-hand side of \eqref{1st-var1} as $\na_v E$ gives the
\emph{first variation formula for the energy}:
\begin{equation} \label{1st-var2}
\na_v E = - \int_M \inn{\tau(\phi), v} \,\om_g \,.
\end{equation}

A smooth map $\phi$ is called \emph{harmonic} if the first variation
\eqref{1st-var2} is zero for all $v \in \Ga(\phi^{-1}TN)$; equivalently,
from \eqref{1st-var2}, 
$\phi$ satisfies the equation
\begin{equation} \label{harmonic}
\tau(\phi) = 0 \,.
\end{equation}
In local coordinates this is an elliptic semilinear system, not linear
except when $(N,h)$ is flat.

Now suppose that $\phi:M \to N$ is harmonic.  Then we can define its
\emph{second variation} as
follows.  By a \emph{(smooth) $2$-parameter variation} $\{\phi_{t,s}\}$
of $\phi$ we mean a smooth map
$M \times I^2 \to N$, \ $(x,t,s) \mapsto \phi_{t,s}(x)$, where
$I^2$ is an open connected subset of $\RR^2$ containing $(0,0)$,
such that $\phi_{0,0} = \phi$.
%For $v,w \in \Ga(\phi^{-1}TN)$.
The \emph{Hessian of
$\phi$} is defined on a pair $v,w$ of vector fields along $\phi$ by
\begin{equation} \label{Hessian1}
H_{\phi}(v,w) =
	\left.\frac{\pa^2 E(\phi_{t,s})}{\pa t\pa s} \right\vert_{(0,0)}
\end{equation}  
where $\{\phi_{t,s}\}$ is a smooth $2$-parameter variation of $\phi$
with
$$
\left.\frac{\pa\phi_{t,s}}{\pa t}\right\vert_{(0,0)} = v \quad \text{and} \quad
\left.\frac{\pa\phi_{t,s}}{\pa s}\right\vert_{(0,0)} = w \,.
$$
Then the left-hand side of \eqref{Hessian1} depends only on $v$ and $w$
and is given by
\begin{equation} \label{Hessian2}
H_{\phi}(v,w) = \int_M \inn{J_{\phi}(v), w} \,\om_g
\end{equation}
where
$$
J_{\phi}(v) = \Delta v - \Tr R^N(\dd\phi, v)\dd\phi \,;
$$
here $\Delta$ denotes the Laplacian on $\phi^{-1}TN$ and $R^N$ the curvature
operator of $N$ (conventions as in \cite{Topics}).
The mapping $J_{\phi}: \Ga(\phi^{-1}TN) \to \Ga(\phi^{-1}TN)$ is called the
\emph{Jacobi operator} (for the energy); it is a self-adjoint
elliptic \emph{linear}
operator.   A vector field $v$ along $\phi$ is called a \emph{Jacobi field (along
$\phi$)} if it is in the kernel of the Jacobi operator, i.e., it
satisfies the \emph{Jacobi equation}
\begin{equation} \label{Jacobi}
J_{\phi}(v) \equiv \Delta v - \Tr R^N(\dd\phi, v)\dd\phi = 0 \,.
\end{equation}
By standard elliptic theory (see, for example, \cite{Mazet}) the set
of Jacobi fields along a
harmonic map is a finite dimensional vector subspace of
$\Ga(\phi^{-1}TN)$.

We shall make use of the following interpretation of the Jacobi operator
as the \emph{linearization} of the tension field.

\begin{proposition} \label{prop:first-order}
Let\/ $\phi:M \to N$ be harmonic and let\/ $v \in \Ga(\phi^{-1}TN)$. 
Let
$\{\phi_t\}$ be a smooth variation of\/ $\phi$ tangent to\/ $v$.  Then
$$
J_{\phi}(v) = - \pat\tau(\phi_t) \,.
$$

In particular, $v$ is a Jacobi field along $\phi$ if and only if the
tension field of $\phi_t$ is zero to first order, i.e.,
\begin{equation} \label{first-order}
\pat\tau(\phi_t) = \tau(\phi) = 0 \,.
\end{equation}
\end{proposition}

\proof
{}From \eqref{Hessian1} we have
\begin{align*}
H_{\phi}(v,w) &= \dt\Bigl( \na_w E(\phi_{t,0}) \Bigr) \\
	&= -\dt \int_M \biginn{\tau(\phi_{t,0}),w} \,\om_g \\
	&= \int_M \Biginn{ - \pat\tau(\phi_{t,0}),w} \,\om_g \,.
\end{align*}
Comparing this with \eqref{Hessian2} gives the statement.
\qed \eproof

We shall write the condition \eqref{first-order} succinctly as
\begin{equation}
\tau(\phi_t) = o(t) \,;
\end{equation}
and we shall call a smooth variation $\{\phi_t\}$ \emph{harmonic to first
order} if it satisfies this condition.  Then the Proposition tells us
that \emph{a smooth variation $\{\phi_t\}$ of a harmonic map $\phi$ is harmonic to
first order if and only if it is tangent to a Jacobi field along $\phi$}.

In particular, if $\{\phi_t\}$ is a smooth variation of $\phi$ through harmonic
maps, its variation vector field
$\ds v = \frac{\pa\phi_t}{\pa t}\Big\vert_0$ is a Jacobi field.
Conversely, a Jacobi field along
a harmonic map is called \emph{integrable} if it is tangent to a
variation through harmonic maps.   Not all Jacobi fields are integrable
as well-known examples along geodesics show; see also \cite{Mu} for two-dimensional
examples. 
Integrability has applications to the study of spaces of harmonic maps and to the
behaviour of a minimizing harmonic map near a singularity,
see  \cite{LemWoo-Jacobi} and the references therein.  It is thus an important
question to determine for what pairs of Riemannian manifolds all Jacobi
fields along harmonic maps are integrable.
The known examples are very few in number; recently L.\ Lemaire and the author have added the case $M = S^2$, \
$N = \CP^2$ \cite{LemWoo-Jacobi}.
One main idea in the proof is that many properties of a
harmonic map are `preserved to first order' under a variation tangent
to a Jacobi field. 
We shall explain this in the next two sections.

\section{Preservation of conformality and isotropy to first order}
We discuss two properties of a harmonic map from a $2$-sphere
to a Riemannian manifold which are preserved to first
order under smooth variations tangent to a Jacobi field; the first,
\emph{conformality}, involves only the first order partial derivatives
of the map, the second, \emph{real isotropy} is a stronger condition
involving higher order partial derivatives.

\subsection{Preservation of conformality}
Let $\phi:(M,g) \to (N,h)$ be a smooth map between Riemannian manifolds.  Then
$\phi$ is called \emph{weakly conformal} if, for each $x \in M$ there is
a number $\la(x) \in [0,\infty)$ such that
$\abs{\dd\phi_x(X)} = \la(x)\abs{X}$ for all $X \in T_xM$.

To formulate this in a way we can use, it is convenient to extend the inner
product on $TN$ to a complex-bilinear inner product $\inn{ \ , \ }^{\cC}$
on the complexified tangent space $T^{\cc}N = TN \otimes_{\RR}\CC$
(and so also on $\phi^{-1}T^{\cc}N = \phi^{-1}TN \otimes_{\RR}\CC$).
We also extend the differential of $\phi$
to a complex-linear map $\dd^{\cc}\phi:T^{\cc}M \to T^{\cc}N$ between
complexified tangent spaces.

Suppose now that $M^2$ is a Riemann surface, i.e. a $1$-dimensional
complex manifold. 
For any complex coordinate $z$ on $M^2$ write
$$
\ds\frac{\pa^{\cc}\phi}{\pa z}
	= (\dd^{\cc}\phi)\Bigl(\frac{\pa}{\pa z}\Bigr) \,,
			\quad
\ds\frac{\pa^{\cc}\phi}{\pa \bar{z}}
	= (\dd^{\cc}\phi)\Bigl(\frac{\pa}{\pa\bar{z}}\Bigr) \,.
$$
Then $\phi$ is weakly conformal if and only if
$$
\eta^{\rR}(\phi) \equiv \Biginn{\frac{\pa^{\cc}\phi}{\pa z},
\frac{\pa^{\cc}\phi}{\pa z}}^{\!\cC} = 0 \,,
$$
this condition being independent of the choice of complex coordinate;
it is equivalent to the vanishing of the $(2,0)$-part of the pull-back $\phi^{-1}h$ of the metric on $N$.

Extend the connection $\DD$ on $\phi^{-1}TN$ by complex-bilinearity
to $\phi^{-1}T^{\cc}N$ and write
$$
\frac{\DD}{\pa z} = \DD_{\ts\frac{\pa}{\pa z}} \,, \quad
\frac{\DD}{\pa\bar{z}} = \DD_{\ts\frac{\pa}{\pa\bar{z}}} \,.
$$
Then, a smooth map $\phi:M^2 \to N$ is harmonic if and only if
\begin{equation} \label{harm2}
\frac{\DD}{\pa\bar{z}}\dcphi = 0 \,, \quad \text{equivalently,} \quad
\frac{\DD}{\pa z}\dcbarphi = 0 \,,
\end{equation}
these expressions giving a non-zero multiple of the tension field
with respect to any Hermitian metric
on $M^2$\,; we thus recover the well-known result that harmonicity is
independent of the choice of such a metric.

Now suppose $\phi:M^2 \to N$ is harmonic.
Then, by the first expression in \eqref{harm2}, 
a smooth variation $\{\phi_t\}$ of $\phi$ is harmonic to first
order if and only if
\begin{equation} \label{ha-1st-order}
\Dt\frac{\DD}{\pa\bar{z}}\dcphit = 0 \,;
\end{equation}
swapping derivatives using the curvature of $N$ shows that
this is equivalent to
$$
\frac{\DD}{\pa\bar{z}}\frac{\DD}{\pa z}\frac{\pa\phi_t}{\pa t}
\Big\vert_{0} +
R\Bigl(\dcbarphi, \frac{\pa\phi_t}{\pa t}\Big\vert_{0} \Bigr) \dcphi
= 0 \,.
$$
It follows that $v \in \Ga(\phi^{-1}TN)$ is a Jacobi field
along $\phi$ if and only if
\begin{equation} \label{Jacobi2}
\frac{\DD}{\pa\bar{z}}\frac{\DD}{\pa z}v
+ R\Bigl(\dcbarphi, v \Bigr)\dcphi = 0 %\,.
\end{equation}
or, equivalently (the conjugate of \eqref{Jacobi2}),

$$
\frac{\DD}{\pa z}\frac{\DD}{\pa\bar{z}}v
+ R\Bigl(\dcphi, v \Bigr)\dcbarphi = 0 \,.
$$
We remark that these equations could be obtained directly from \eqref{Jacobi}.
Indeed, their sum is that equation
and their difference is a consequence of the
first Bianchi identity for the curvature.

Now let $\phi:M^2 \to N$ be a smooth map and let $\{\phi_t\}$ be a
smooth variation of it. 
For any complex coordinate, set
$$
\eta^{\rR}(\phi_t) = \Biginn{\frac{\pa^{\cc}\phi_t}{\pa z},
	\frac{\pa^{\cc}\phi_t}{\pa z}}^{\!\cC}.
$$
Say that $\{\phi_t\}$ is \emph{conformal to first order} if
$\eta^{\rR}(\phi_t) = o(t)$,
i.e.
$$
\pat\eta^{\rR}(\phi_t) = \eta^{\rR}(\phi) = 0 \,.
$$   

Let $v = \ds\left.\frac{\pa\phi_t}{\pa t} \right\vert_{0}$ be the variation
vector field of $\{\phi_t\}$.  Then, since $\DD$ is a Riemannian
connection, 
\begin{equation} \label{conf}
\pat\eta^{\rR}(\phi_t) =
   2\Biginn{\frac{\DD v}{\pa z}, \frac{\pa^{\cc}\phi}{\pa z}}^{\!\cC}.
\end{equation}
We call a vector field $v$ along a smooth map $\phi:M^2 \to N$
\emph{conformal} if
$\ds\Biginn{\frac{\DD v}{\pa z}, \frac{\pa^{\cc}\phi}{\pa z}}^{\!\cC}
= 0$; this condition is clearly independent of the choice of complex
coordinate.  Then from \eqref{conf} we see that conformality to first
order is equivalent to conformality of the variation vector field
as follows:

\begin{lemma} \label{lem:WC}
Let\/ $\phi:M^2 \to N$ be a weakly conformal map from a Riemann surface to a
Riemannian manifold and let\/ $v \in
\Ga(\phi^{-1}TN)$.  Then the following are equivalent:

\begin{enumerate} \vspace{-1mm}
\item[{\rm (i)}] $v$ is a conformal vector field along\/ $\phi$; \vspace{-2mm}
\item[{\rm (ii)}] \emph{there exists} a smooth variation of\/ $\phi$
tangent to\/ $v$ which is conformal to first order; \vspace{-2mm}
\item[{\rm (iii)}] \emph{all} smooth variations of\/ $\phi$
tangent to\/ $v$ are conformal to first order.
\end{enumerate}
\vspace{-7mm} \qed \end{lemma}

Any harmonic map from the $2$-sphere $S^2$ to an arbitrary
Riemannian manifold $N$ is weakly conformal.  This is proved by
showing that, for any harmonic map $\phi:M^2 \to N$ from a Riemann
surface, $\eta^{\rR}(\phi)\,\dd z^2$ defines a
\emph{holomorphic quadratic differential}, i.e. a
holomorphic section of the holomorphic line bundle $\otimes^2 T'_*M^2$;
if $M^2 = S^2$, this bundle
has negative degree and so any holomorphic section vanishes, see,
for example, \cite{Wood-Trieste}.  We
generalize this in the next result; for later comparison we give two
proofs. 

\begin{proposition} \label{prop:conf}
Any Jacobi field along a harmonic map $\phi:S^2 \to N$ is conformal;
equivalently, any smooth variation of $\phi$ tangent to a Jacobi field
is conformal to first order.
\end{proposition}

\proof \!\!{\bf 1} \
Let $\phi:M^2 \to N$ be a harmonic map from a Riemann surface and let
$v$ be a Jacobi field along it.
Using the Jacobi equation \eqref{Jacobi2}, in any complex coordinate we have
\begin{align*}
\dbar\Biginn{\frac{\DD v}{\pa z},\frac{\pa^{\cc}\phi}{\pa z}}^{\!\cC}
&= \Biginn{\frac{\DD}{\pa\bar{z}}\frac{\DD v}{\pa z},
		\frac{\pa^{\cc}\phi}{\pa z}}^{\!\cC}
	+ \Biginn{\frac{\DD v}{\pa z},
		\frac{\DD}{\pa\bar{z}}\frac{\pa^{\cc}\phi}{\pa z}}^{\!\cC}
\\ &= -  \Biginn{R\Bigl(\frac{\pa^{\cc}\phi}{\pa\bar{z}}, v \Bigr)
	\frac{\pa^{\cc}\phi}{\pa z},
\frac{\pa^{\cc}\phi}{\pa z}}^{\!\cC} + 0 
\\ &=0 \,.
\end{align*}
Hence
$\ds\Biginn{\frac{\DD v}{\pa z},\frac{\pa^{\cc}\phi}{\pa z}}^{\!\cC}\dd z^2$
defines a holomorphic differential on $M^2$; if $M^2 = S^2$ this must
vanish.

\medskip

\proof \!\!{\bf 2} \
Let $\phi:M^2 \to N$ be a harmonic map from a Riemann surface and let
$\{\phi_t\}$ be a smooth variation of it tangent to a Jacobi field.  Then, by Proposition
\ref{prop:first-order},
$$
\dbar\eta^{\rR}(\phi_t) =
	2\Biginn{\frac{\DD}{\pa\bar{z}}\frac{\pa^{\cc}\phi_t}{\pa z},
		\frac{\pa^{\cc}\phi_t}{\pa z}}^{\!\cC} = o(t) \,.
$$
It follows that
$$
\dbar\pat\eta^{\rR}(\phi_t) =
	\pat\dbar\eta^{\rR}(\phi_t) = 0 \,.
$$
Hence, as well as $\eta^{\rR}(\phi)$,
$\ds\pat\eta^{\rR}(\phi_t)$ is a holomorphic
quadratic differential and so must vanish if $M^2 = S^2$.
\qed \eproof

\subsection{Preservation of real isotropy}
For any smooth map $\phi:M^2 \to N$ from a Riemann surface to a Riemannian
manifold, integers $r,s \geq 1$ and local complex coordinate $z$, set
$$
\eta^{\rR}_{r,s}(\phi) = \Biginn{\Dz[r-1]\dcphi, \Dz[s-1]\dcphi}^{\!\cC}.
$$
Recall that $\phi$ is called \emph{real isotropic} or
\emph{pseudoholomorphic} \cite{Ca2} if
$\eta^{\rR}_{r,s}(\phi) = 0$ for all integers $r,s \geq 1$.  Note
that this condition is
independent of the complex coordinate chosen and that any real isotropic
map is weakly conformal.

Let $\phi:S^2 \to N$ be a harmonic map from the $2$-sphere to a
\emph{space form}, i.e.\ a Riemannian manifold of
constant sectional curvature.    Then \cite{Ca1,Ca2}
$\phi$ is real isotropic.  This is proved by
showing, inductively that, for $k =2,3,\ldots$ and $r+s=k$,
$\eta^{\rR}_{r,s}(\phi)\dd z^k$ defines a
\emph{holomorphic $k$-differential}, i.e., a holomorphic section of
$\otimes^kT'_*S^2$, and so
must vanish.   It follows that all harmonic maps can be constructed explicitly
from holomorphic maps, see \cite{Ca2,Chern, Barbosa,Bor-Gar,Lawson}.
We now generalize this argument.  Let $\{\phi_t\}$ be a
$1$-parameter variation of $\phi$.   Say that \emph{$\{\phi_t\}$ is
real isotropic to first order} if $\eta^{\rR}_{r,s}(\phi_t) = o(t)$ for
all $r,s$, i.e.,
$$
\pat \eta^{\rR}_{r,s}(\phi_t) = \eta^{\rR}_{r,s}(\phi) = 0 \qquad (r,s \geq 1) \,.
$$
Note that this condition is independent of the choice of complex
coordinate. It can be formulated in terms of the variation vector field
of $\{\phi_t\}$ as follows.   For any smooth map $\phi:M^2 \to N$ from a
Riemann surface to a Riemannian manifold, vector field $v$ along $\phi$,
integers $r,s \geq 1$ and local complex coordinate $z$, set
$$
j^{\rR}_{r,s}(v)
	= \Biginn{\frac{\DD^r v}{\pa z^r}, \Dz[s-1]\dcphi}^{\!\cC}
	+ \Biginn{\Dz[r-1]\dcphi, \frac{\DD^s v}{\pa z^s}}^{\!\cC} \,.
$$
Then we have

\begin{lemma} \label{lem:j-real}
Let\/ $\phi:M^2 \to N$ be a real isotropic map from a Riemann surface into a
space form and let\/ $v$ be a Jacobi field along it.   Then, for any
smooth variation\/ $\{\phi_t\}$ of\/ $\phi$ tangent to\/ $v$
and any complex coordinate\/ $z$ on\/ $M^2$ we have
$$
j^{\rR}_{r,s}(v) = \pat\eta^{\rR}_{r,s}(\phi_t) \qquad (r,s \geq 1) \,.
$$
\end{lemma}

\proof
We show by induction that, for $k = 1,2, \ldots$,
\begin{equation} \label{commute}
\Dt\Dz[k-1]\dcphit \in \theta_{k-1}(\phi) + \frac{\DD^k v}{\pa z^k}
\end{equation}
where
$$
\theta_{k-1}(\phi) = \text{span}\left\{ \frac{\pa^{\cc}\phi}{\pa z},
\ldots, \Dz[k-2]\dcphi \right\} \ (k \geq 2) \  \text{ and }
\ \theta_{0}(\phi) = \{\vec{0}\} \,;
$$
note that this is independent of the choice of complex coordinate.

For $k=1$, \eqref{commute} holds since
$\ds\Dt \frac{\pa^{\cc}\phi_t}{\pa z}
	= \left.\frac{\DD}{\pa z}\frac{\pa\phi_t}{\pa t}\right\vert_{0}
	= \frac{\DD v}{\pa z}$\,.

Suppose that \eqref{commute} holds for $k-1$ for some $k \geq 2$.   Then
\begin{equation} \label{induct}
\Dt \Dz[k-1]\dcphit = \frac{\DD}{\pa z}\Dt\Dz[k-2]\dcphit + R\Bigl(
\frac{\pa^{\cc}\phi}{\pa z}, v \Bigr)\Dz[k-2]\dcphi \,.
\end{equation}
Now, by the induction hypothesis the first term on the right-hand side
of \eqref{induct} lies in
$\ds \theta_{k-1}(\phi) + \frac{\DD^k v}{\pa z^k}$\,. 
Also, by the well-known formula for the curvature of a
space form (see, for example, \cite[Chapter V]{K-N1}), the curvature
term is a multiple of
\begin{equation} \label{curv0}
\Biginn{v, \Dz[k-2]\dcphi}^{\!\cC} \, \frac{\pa^{\cc}\phi}{\pa z}
	- \Biginn{\frac{\pa^{\cc}\phi}{\pa z},
\Dz[k-2]\dcphi}^{\!\cC}v \,;
\end{equation}
since the second term of this is zero by the isotropy of $\phi$, 
\eqref{curv0} lies in $\theta_{1}(\phi)$ and so in
$\theta_{k-1}(\phi)$.  It follows that the right-hand side of
\eqref{induct} lies in  
$\ds\theta_{k-1}(\phi) + \frac{\DD^k v}{\pa z^k}$ as required, completing
the inductive step.

Using \eqref{commute} and the isotropy of $\phi$ we have, for any $r,s \geq 1$,
\begin{align*}
\pat\eta^{\rR}_{r,s}(\phi_t) 
	&= \Biginn{\Dt\frac{\DD^{r-1}}{\pa z^{r-1}}
	\frac{\pa^{\cc}\phi_t}{\pa z},
	\Dz[s-1]\dcphi}^{\!\cC}
	+ \Biginn{\Dz[r-1]\dcphi,
\Dt\frac{\DD^{s-1}}{\pa z^{s-1}}\frac{\pa^{\cc}\phi_t}{\pa z}}^{\!\cC}\\
	&= \Biginn{\frac{\DD^r v}{\pa z^r}, \Dz[s-1]\dcphi}^{\!\cC}
	+ \Biginn{\Dz[r-1]\dcphi, \frac{\DD^s v}{\pa z^s}}^{\!\cC} \\
	&= j^{\rR}_{r,s}(v)
\end{align*}
as required.
\qed \eproof

Say that a vector field $v$ along a smooth map $\phi:M^2 \to N$
\emph{preserves real isotropy to first order} if $j^{\rR}_{r,s}(v) = 0$ for
all $r,s \geq 1$; note that this condition is independent of the
choice of local complex coordinate.

\begin{proposition} \label{prop:pres-real}
Let $\phi:S^2 \to N$ be a harmonic map from the $2$-sphere to a space
form.  Then any Jacobi field preserves real isotropy to first order,
equivalently, any smooth variation of\/ $\phi$ tangent to a
Jacobi field is real isotropic to first order.
\end{proposition}

\proof
Let $\{\phi_t\}$ be any smooth variation of $\phi$.   We shall show that
$\eta^{\rR}_{r,s}(\phi_t) = o(t)$ for all integers $r,s \geq 1$.
In fact we show by induction that, for all $K \in \{1,2,\ldots\}$,
\begin{align}
{\rm (i)} \quad &\frac{\DD}{\pa\bar{z}}\Dz[K-1]\dcphit \in
	\theta_{K-1}(\phi_t) + o(t) \,, \label{ind1} \\ 
{\rm (ii)} \quad &\text{for all } r,s \geq 1 \text{ with } r+s = K+1, \
	\eta^{\rR}_{r,s}(\phi_t) = o(t) \label{ind2} \,.
\end{align}
This is true for $K=1$ by \eqref{ha-1st-order} and Proposition \ref{prop:conf}.

Suppose that it is true for all $K < k$ for some $k \geq 2$.  We shall show
that it is true for $K = k$.   Indeed,
\begin{equation} \label{curv}
\frac{\DD}{\pa\bar{z}}\Dz[k-1]\dcphit =
R\Bigl(\frac{\pa^{\cc}\phi_t}{\pa z},
\frac{\pa^{\cc}\phi_t}{\pa\bar{z}}\Bigr)\Dz[k-2]\dcphit
	+ \frac{\DD}{\pa z}\frac{\DD}{\pa\bar{z}}\Dz[k-2]\dcphit \,.
\end{equation}
The curvature term is a multiple of
$$
\Biginn{\frac{\pa^{\cc}\phi_t}{\pa\bar{z}},
	\Dz[k-2]\dcphit}^{\!\cC} \, \frac{\pa^{\cc}\phi_t}{\pa z} 
	- \Biginn{\frac{\pa^{\cc}\phi_t}{\pa z},
	\Dz[k-2]\dcphit}^{\!\cC} \, \frac{\pa^{\cc}\phi_t}{\pa\bar{z}} \,.
$$
The first term of this lies in $\theta_{1}(\phi_t)$ and so in
$\theta_{k-1}(\phi_t)$; by the induction hypothesis, the second term is
$o(t)$.   Also, by the induction hypothesis, the last term of \eqref{curv}
lies in $\theta_{k-1}(\phi_t)$, so that \eqref{ind1} holds for $K=k$.
Next, for $r,s \geq 1$ with $r+s = k+1$ we have
$$
\frac{\pa}{\pa\bar{z}}\eta^{\rR}_{r,s}(\phi_t) =
 \Biginn{\frac{\DD}{\pa\bar{z}}\Dz[r-1]\dcphit,\Dz[s-1]\dcphit}^{\!\cC}
 + \Biginn{\Dz[r-1]\dcphit,\frac{\DD}{\pa\bar{z}}\Dz[s-1]\dcphit}^{\!\cC}
\,.
$$
By the induction hypothesis, this is $o(1)$, hence
$\ds\pat\eta^{\rR}_{r,s}(\phi_t)\dd z^{k+1}$ defines a holomorphic
$(k+1)$-differential on $S^2$.
As before, this must vanish, so that $\eta^{\rR}_{r,s}(\phi_t) =
o(t)$ completing the induction step.
\qed \eproof

\section{Preservation of holomorphicity and complex isotropy
to first order}
We give two properties of a harmonic map into a K\"ahler
manifold which are preserved
to first order under smooth variations tangent to a Jacobi field;
the first is holomorphicity which involves only first order partial
derivatives, the second is complex isotropy which involves higher order
partial derivatives.

\subsection{Preservation of holomorphicity}
Let $M$ and $N$ be compact K\"ahler manifolds.
Then we can decompose the complexified tangent spaces $T^{\cc}M = TM
\otimes_{\RR}\CC$ and $T^{\cc}N = TN \otimes_{\RR}\CC$ into $(1,0)$ and
$(0,1)$ parts, viz.,
\begin{equation} \label{decomp}
T^{\cc}M = T'\!M \oplus T''\!M \qquad T^{\cc}N = T'\!N \oplus T''\!N \,.
\end{equation}
We shall let $\inn{ \ , \ }^{\hh}$ denote the \emph{Hermitian extension} of the
inner product on $TN$ to $T^{\cc}N$ given by
$$
\inn{v , w}^{\hh} = \inn{v , \ov{w}}^{\cC} \,;
$$
this restricts to a positive definite Hermitian inner product on $T'\!N$
and on its pull-back $\phi^{-1}T'\!N$.
Let $\phi:M \to N$ be a smooth map.  As before extend its
differential to a complex-linear map $\dd^{\cc}\phi:T^{\cc}M \to
T^{\cc}N$; then we can consider the two components:
$$
\pa\phi:T'\!M \to T'\!N \,, \qquad \bar{\pa}\phi:T''\!M \to T'\!N \,.
$$
Recall that $\phi$ is called \emph{holomorphic} if $\bar{\pa}\phi = 0$
(and \emph{antiholomorphic} if $\pa\phi = 0$).

For simplicity we calculate in a local complex coordinate system
$(z^j)$ on $M$, writing
$$
\frac{\pa\phi}{\pa z^j} =
	(\pa\phi)\bigl(\frac{\pa}{\pa z^j}\bigr) \,,
		\quad
\frac{\pa\phi}{\pa z^{\ov{j}}} =
	(\bar{\pa}\phi)\bigl(\frac{\pa}{\pa z^{\ov{j}}}\bigr) \,;
$$
note that
$$
\frac{\pa^{\cc}\phi}{\pa z^j}
	= \frac{\pa\phi}{\pa z^j} + \ov{\frac{\pa\phi}{\pa z^{\ov{j}}}} \,.
$$
Say that a smooth variation $\{\phi_t\}$ is \emph{holomorphic to
first order} if $\bar{\pa}\phi_t = o(t)$, i.e., for any local complex
coordinate system $(z^j)$\,,
$$
\left. \frac{\na}{\pa t}\right\vert_{0}
	\!\!\frac{\pa\phi_t}{\pa z^{\ov{j}}}
	= \frac{\pa\phi}{\pa z^{\ov{j}}} = 0 \quad \forall j \,.
$$

For any $v \in \Ga(\phi^{-1}T^{\cc}N)$, let $v'$
denote its $(1,0)$ component under the decomposition \eqref{decomp}.
Note that, if $v$ is real, i.e., lies in $\Ga(\phi^{-1}TN)$,
$v'$ determines $v$, indeed $v = v' + \ov{v'}$.  
Say that $v$ is \emph{holomorphic} if $\DD_Z v' = 0$ for all
$Z \in T''M$, equivalently,
\begin{equation}
\frac{\DD v'}{\pa z^{\ov{j}}} = 0 \quad \forall j \,.
\end{equation}

Holomorphicity to first order is then equivalent to holomorphicity of
the variation vector field as follows (cf.\ Lemma \ref{lem:WC}):

\begin{lemma} \label{lem:hol}
Let\/ $\phi:M \to N$ be a holomorphic map between K\"ahler manifolds and let\/
$v \in \Ga(\phi^{-1}TN)$.  Then the following are equivalent:

\begin{enumerate} \vspace{-1mm}
\item[{\rm (i)}] $v$ is a holomorphic vector field along\/ $\phi$; \vspace{-2 mm}
\item[{\rm (ii)}] \emph{there exists} a smooth variation of\/ $\phi$
tangent to $v$ which is holomorphic to first order; \vspace{-2 mm}
\item[{\rm (iii)}] \emph{all} smooth variations of\/ $\phi$
tangent to $v$ are holomorphic to first order.
\end{enumerate}
\vspace{-7mm} \qed \end{lemma}

\proof
Noting that the K\"ahler condition on $N$ implies that $(\DD_Z
v)' = \DD_Z(v')$ for all $v \in \Ga(\phi^{-1}T^{\cc}N)$, \ $Z \in
T^{\cc}M$,
we have
$$
\Dt \frac{\pa\phi_t}{\pa z^{\ov{j}}} =
\frac{\DD}{\pa z^{\ov{j}}} \left. \left(\frac{\pa\phi_t}{\pa
t}\right\vert_{0} \right)' = \frac{\DD v'}{\pa z^{\ov{j}}} \quad \forall j
$$
from which the lemma follows.
\qed \eproof

\begin{proposition} \label{prop:hol}
Let\/ $\phi:M \to N$ be a holomorphic map between compact K\"ahler
manifolds.
Then any Jacobi field along\/ $\phi$ is holomorphic; equivalently, any
smooth variation of\/ $\phi$ tangent to a Jacobi field is holomorphic to
first order.
\end{proposition}

\proof
The harmonicity condition for a smooth map $\phi:M \to N$ between
K\"ahler manifolds is (with summation over $i$ and $j$) \ 
$\ds g^{i\ov{j}} \frac{\na}{\pa z^i} \frac{\pa\phi}{\pa z^{\ov{j}}}
	= 0$\,.
Hence a smooth variation $\{\phi_t\}$ of a harmonic map $\phi:M \to N$
is tangent to a Jacobi field $v$ if and only if
$$
\pat g^{i\ov{j}} \frac{\na}{\pa z^i} \frac{\pa\phi_t}{\pa z^{\ov{j}}}
	= 0 \,.
$$
On using the curvature to commute the derivatives, this
gives the Jacobi equation in the form
$$
g^{i\ov{j}} \left\{ \frac{\na^2 v'}{\pa z^i \pa z^{\ov{j}}}
	+ R\Bigl( \frac{\pa^c\phi}{\pa z^i}, v\Bigr)
\frac{\pa\phi}{\pa z^{\ov{j}}} \right\} = 0 \,.
$$

If $\phi$ is holomorphic, $\pa\phi/\pa z^{\bar j} = 0$ so that
the curvature term vanishes; then integrating by parts gives
$$
\int_M g^{i\ov{j}}
\Biginn{\frac{\DD v'}{\pa z^{\ov{j}}},\frac{\DD v'}
	{\pa z^{\ov{i}}}}^{\!\hh} \,\om_g = 0 \,.
$$
Therefore $\ds\frac{\DD v'}{\pa z^{\ov{j}}} = 0$ \ $\forall j$,
so that $v$ is holomorphic.
\qed \eproof

For an alternative proof see \cite{LemWoo-Jacobi}.

\subsection{Preservation of complex isotropy}
Let $N$ be a K\"ahler manifold.  For any smooth map $\phi:M^2
\to N$ from a Riemann surface, integers $r,s \geq 1$ and complex
coordinate $z$, set
$$
\eta^{\cC}_{r,s}(\phi)
	= \Biginn{\Dz[r-1]\dphi, \Dzbar[s-1]\dbarphi}^{\!\hh} \,.
$$
Then $\phi$ is called \emph{complex isotropic} if
$\eta^{\cC}_{r,s}(\phi) = 0$ for all integers $r,s \geq 1$.  Note that
this condition is
independent of the complex coordinate chosen and that it is stronger
than real isotropy.

It is well-known that any harmonic map $\phi:S^2 \to N$ from the $2$-sphere to a complex space
form (i.e.\ a K\"ahler manifold of constant holomorphic sectional curvature) is
complex isotropic, see, for example, \cite{D-Z,E-W,Wood-Horwood}.
As in the real case, this is proved by
showing inductively that, for $k =2,3,\ldots$ and $r+s=k$,
$\eta^{\cC}_{r,s}(\phi)\dd z^k$
defines a holomorphic $k$-differential on $S^2$ and so must vanish.

For any smooth map $\phi:M^2 \to N$ from a
Riemann surface to a K\"ahler manifold, vector field $v$ along $\phi$,
integers $r,s \geq 1$ and local complex coordinate $z$, set
$$
j^{\cC}_{r,s}(v)
	= \Biginn{\frac{\DD^r v}{\pa z^r}, \Dzbar[s-1]\dbarphi}^{\!\hh}
	+ \Biginn{\Dz[r-1]\dphi, \frac{\DD^s v}{\pa\bar{z}^s}}^{\!\hh} \,.
$$
Then, analogously to Lemma \ref{lem:j-real} we have

\begin{lemma} 
Let\/ $\phi:M^2 \to N$ be a complex isotropic map from a Riemann surface into a
space form and let\/ $v$ be a Jacobi field along it.   Then, for any
smooth variation\/ $\{\phi_t\}$ of\/ $\phi$ tangent to\/ $v$
and any complex coordinate\/ $z$ on\/ $M^2$,
$$
j^{\cC}_{r,s}(v) = \pat\eta^{\cC}_{r,s}(\phi_t) \qquad (r,s \geq 1) \,.
$$
\qed \end{lemma}

Say that \emph{$v$ preserves complex isotropy to first order}
if $j^{\cC}_{r,s}(v) = 0$ for
all $r,s \geq 1$, and say that \emph{$\{\phi_t\}$ is
complex isotropic to first order} if
$\eta^{\cC}_{r,s}(\phi_t) = o(t)$, i.e.,
$$
\pat \eta^{\cC}_{r,s}(\phi_t) = \eta^{\cC}_{r,s}(\phi) = 0
	\qquad (r,s \geq 1) \,.
$$
Note that these conditions are independent of the choice of local
complex coordinate.   Then we have the complex analogue of Proposition
\ref{prop:pres-real}:

\begin{proposition} \label{prop:pres-cx}
Let\/ $\phi:S^2 \to N$ be a harmonic map from the $2$-sphere to a complex space
form.  Then any Jacobi field preserves complex isotropy to first order,
equivalently, any smooth variation of\/ $\phi$ tangent to a Jacobi
field is complex isotropic to first order.
\qed \end{proposition}

The proofs proceed in the same way as for the real case, except that we
must use the more complicated formula for the curvature of a complex
space form (see, for example, \cite[Chapter IX]{K-N2}).  The
calculations are omitted; they are simplified by the
following easy consequence of that curvature formula:

\begin{lemma}
For any $X \in T^{\cC}N$, \ $\ds R\!\!\left(X, \dcphi\right)\Dz[k-1]\dphi$ \ is a
linear combination of\/ \,$\ds \dphi$\,,\, $\ds \Dz[k-1]\dphi$\, and\/\,
$\ds \eta^{\cC}_{k,1}(\phi)X'$.
\qed \end{lemma}

\end{document}